\documentclass{article}
\usepackage{amsmath}
\usepackage{amssymb}
\usepackage{amsthm}
\usepackage{graphicx}

\newtheorem{thm}{Theorem}

\newtheorem{lemma}[thm]{Lemma}

\title{Convergence of the Freely Rotating Chain to the Kratky-Porod Model of Semi-flexible Polymers}

\author{Philip Kilanowski,\thanks{Dominican House of Studies, Washington DC}
\hspace{0.3cm} Peter March,\thanks{Rutgers University, Piscataway, NJ}$\hspace{.2cm}^{\S}$
\hspace{0.3cm} Marko \v{S}amara \thanks{Arizona State University, Tempe AZ \newline 
\indent \hspace{.02cm} $\S$\hspace{.02cm}This research was supported by an Intergovernmental Personnel Act Mobility Assignment NSF DMS 0649591}}

\begin{document}

\maketitle

\begin{abstract}
\noindent The freely rotating chain is one of the classic discrete models of a polymer in dilute solution. It consists of a broken line of $N$ straight segments of fixed length such that the angle between adjacent segments is constant and the $N-1$ torsional angles are independent, identically distributed, uniform random variables. We provide a rigorous proof of a folklore result in the chemical physics literature stating that under an appropriate scaling, as $N\rightarrow\infty$, the freely rotating chain converges to a random curve defined by the property that its derivative with respect to arclength is a brownian motion on the unit sphere. This is the Kratky-Porod model of semi-flexible polymers. We also investigate limits of the model when a stiffness parameter, called the persistence length, tends to zero or infinity. The main idea is to introduce orthogonal frames adapted to the polymer and to express conformational changes in the polymer in terms of stochastic equations for the rotation of these frames.
\end{abstract}



\pagestyle{myheadings}


\section{Introduction}
The simplest caricature of a polymer in dilute solution is a sequence of random vectors $R_0, R_1, \dots, R_N$ in $\mathbb R^3$ joined together by line segments $Q_n=R_n-R_{n-1},~1\leq n\leq N$. The freely jointed chain, or bead/rod model, introduced by Kramers in 1946 \cite{Kr46}, assumes the segments, or bonds, $Q_n$ are independent, identically distributed, random variables uniformly distributed on a sphere of radius $a$. A closely related variant, introduced by Rouse in 1953 \cite{Ro53}, called the bead/spring model, assumes the segments $Q_n$ are i.i.d gaussian random vectors. The physical interpretation in Kramers' model is that statistically significant sections of the polymer behave as if they were inextensible rods that can rotate freely with respect to one another. In Rouse's model, the physical interpretation is that statistically significant sections of the polymer behave as if they were linear springs connecting $N+1$ beads in which the thermal forces of the solvent acting on the beads are equilibrated by the elastic restoring forces of the springs.\\

\noindent These models provide a reasonable description of ideal or flexible polymers, in which correlations between monomers decay rapidly with distance along the polymer. That they provide equivalent descriptions of the large scale properties of flexible polymers can be seen from the fact that the bead/rod and bead/spring models are $N$ step random walks that converge to three dimensional brownian motion, as $N\rightarrow\infty$, under a suitable scaling. However, these are poor models for stiff or semi-flexible polymers, in which the monomer-monomer correlation length is comparable to the length of the polymer itself. To address this issue, Kirkwood and Riseman \cite{KR48}, \cite{KR49} introduced in 1948 a variant of the Kramers chain, called the \emph{freely rotating chain}, in which geometric constraints are imposed on the segments so that they are correlated with one another and no longer independent. In this model, the bonds have a common length $a$; the angle between adjacent bonds, or bond angle, has a fixed value $\theta$; and the bond $Q_n$ is chosen independently and uniformly at random from the cone of axis $Q_{n-1}$, aperture $\theta$ and side length $a$, (cf. Figure 1). While the bonds $Q_n$ are no longer independent, they do form a Markov chain, namely, a geodesic random walk on the sphere of radius $a$ of step size $a\theta$. We recommend standard references to the classical polymer physics literature such as \cite{Flory}, \cite{DE}, \cite{Volkenstein} and \cite{Yamakawa} for more in depth discussion of the basic physical phenomena and the standard mathematical models of polymers.

\begin{figure}[h!]
\label{conformation}
\begin{center}
\includegraphics[scale=0.65]{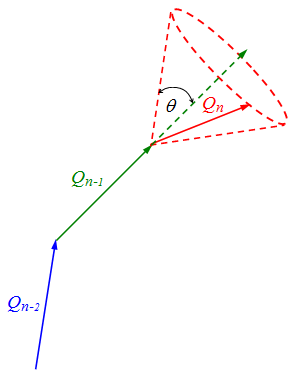}\label{segments}
\caption{Construction of the discrete polymer from its bonds $Q_n$}
\end{center}
\end{figure}

\noindent At the same time, a continuum version of the freely rotating chain was proposed independently by Kratky and Porod \cite{Kratky and Porod}. Their model is defined by a Hamiltonian for a curve $R(s)$ in terms of its derivative with respect to arclength $u(s)=\frac{\partial R(s)}{\partial s},$
$$
H(R) = \lambda\int_0^L\left|\frac{\partial u(s)}{\partial s}\right|^2ds.
$$
Here $\lambda$ is a physical constant, and $L$ is the arclength of the polymer. By definition, the tangent vector $u(s)$ has unit length, but it is often pointed out explicitly that $|u(s)|\equiv 1$, because in the physical formalism, this fact greatly complicates the implied functional integration over paths that is necessary to obtain the partition function of the model. For this reason, the Kratky-Porod model is considered to be somewhat difficult to work with.\\

\noindent A formal calculation relates the Kratky-Porod model to the scaling limit of the freely rotating model. Since

$$
|Q_k-Q_{k-1}|^2=2a^2(1-\cos\theta),
$$\

\noindent it follows: first, that the difference quotients $(Q_k-Q_{k-1})/a$ approximate a derivative with respect to arclength, provided $a=L/N$; and second, that the sum of squares of these difference quotients is bounded if $\theta=O(1/\sqrt{N})$; say, $\theta=\kappa/\sqrt{N}$. This formal calculation is the basis of a folklore result in the chemical physics literature that the Kratky-Porod model is the limit of the freely jointed model under the indicated scaling. Now, it is difficult to make rigorous sense of the Hamiltonian formalism of Kratky and Porod, and we do not try to do so. Rather, we show that under this scaling, the freely jointed chain converges to a random curve whose derivative with respect to arclength can be identified explicitly as a spherical brownian motion. This is easy to guess, referring back to the probabilistic description of $Q_n$ as a geodesic random walk on the sphere of radius $a$. Evidently, the normalized segments $Q_n/a$ lie on the unit sphere and the steps $Q_{n-1}\rightarrow Q_{n}$ have uniformly distributed directions and are taken along great circles of arclength $O(1/\sqrt{N})$, which is the correct scaling for convergence to spherical brownian motion. The first result of this article makes these arguments precise and rigorously defines the Kratky-Porod model.\\

\noindent Note that throughout this work, we use pinned boundary conditions at $s=0$. Thus, the initial bead is always $R_0=0$, and the initial bond is always in the vertical direction in the laboratory frame, namely, $Q_1=ae_3$ where $e_3=(0,0,1).$

\begin{thm}
\label{Theorem 1}
Let $R^N_n,~0\leq n\leq N$, be the freely rotating chain with bond length $a=L/N$ and bond angle $\theta=\kappa/\sqrt{N}$. Let $R^N(s),~ 0\leq s\leq L$ be the piecewise linear curve obtained by linear interpolation of the beads $R^N(nL/N)=R^N_n.$ Let $Q_s(\lambda),~0\leq s\leq L$ be brownian motion on the unit sphere $S^2$, generated by $\frac{1}{\lambda}\Delta_{S^2}$, starting from the north pole $Q_0(\lambda)=e_3=(0,0,1)$ and let $\ell_p=2L/\kappa^2$. Then the processes $R^N$ converge weakly, as $N\rightarrow\infty$, to the Kratky-Porod model defined by

$$
R_s=\int_0^sQ_u(\ell_p)du,~ 0\leq s\leq L.
$$
\end{thm}

\noindent Note that the unit tangent process $Q_s(\ell_p)$ is not standard brownian motion on $S^2$ since it is generated by $\frac{1}{\ell_p}\Delta_{S^2}$ rather than the standard $\frac{1}{2}\Delta_{S^2}$. The reason for the non-standard normalization is to focus attention on the parameter $\ell_p=2L/\kappa^2$, called the persistence length. Since

\begin{align}
E[Q_s(\ell_p)\cdot Q_t(\ell_p)]=e^{-2|t-s|/\ell_p}
\end{align}
as a standard calculation shows (see Section 2), $\ell_p$ is a correlation length that describes the exponential decay of tangent-tangent correlations along the polymer. Values of $\ell_p$ that are large compared to the polymer length $L$ are indicative of stiff, rod-like polymers, whereas values $\ell_p$ that are small compared to $L$ are indicative of flexible, random coil polymers. Thus, the dimensionless ratio

$$
\frac{L}{\ell_p}=\frac{\kappa^2}{2}
$$
characterizes the intrinsic stiffness of the polymer and gives a physical interpretation of the bond angle parameter $\kappa$. The persistence length is the key parameter describing a qualitative transition in the Kratky-Porod model which is recorded in the second main result of this article.

\begin{thm}
\label{Theorem 2}
The Kratky-Porod model with fixed length $L$ exhibits the hard rod/random coil transition, in the sense that:

A. (Hard Rod) As $\ell_p\to\infty$, $(R_s; \, 0\le s \le L)\, \to \, (se_3;\, 0\le s \le L)$ in probability, and $(\sqrt{\ell_p}(R_s-se_3);\, 0\le s \le L)$ converges weakly to a Gaussian process $(W_s;\, 0 \le s \le L)$ with components $W^3_s\equiv 0$, and for $i=1,2$,

$$
W^i_s=\int_0^s\beta^i_udu,~ 0\leq s\leq L  
$$

\noindent where the $\beta^i$'s are independent standard one dimensional brownian motions starting at zero.

B. (Random Coil) As $\ell_p\to 0$,  $(\sqrt{3/\ell_p} R_s;\, 0\le s \le L)$ converges weakly to $(W_s;\, 0 \le s \le L)$, where $W$ is a standard 3-dimensional Brownian motion.

\end{thm}

\noindent The significance of Theorem 2 is that it shows how the Kratky-Porod model interpolates between two extreme polymer types. It confirms, as expected, that the freely rotating chain model for semi-flexible polymers agrees with the freely jointed chain and the bead/spring chain in the large $N$ and small $\ell_p$ regime. It also confirms the expected hard rod limit as $\ell_p\to\infty$, but the precise form of the gaussian fluctuations on the scale $O(\sqrt{\ell_p})$ seems to be new. This regime is known as the weakly bending limit in the polymer physics literature, e.g., \cite{DE}.\\

\noindent The main idea in the proofs of these results is a representation $Q_n=Z_nQ_1$ of the $n^{th}$ bond of the freely rotating chain as a random rotation of the initial bond. The Markovian structure of the bond process $Q_n$ leads to a simple representation of the rotation process $Z_n$ in terms of i.i.d. rotations. It develops the fact that $Z_n$ is the solution of a simple stochastic difference equation in the Lie group $SO(3)$ that corresponds, in the scaling limit, to the Stratonovich SDE,

\begin{align}
\partial Z_s=\frac{1}{\sqrt{\ell_p}}Z_s\partial B_s
\end{align}
Here $B(s)$ is a brownian motion in the Lie algebra of 3-by-3 anti-symmetric matrices given by

\begin{equation}
B(s)=\left(\begin{array}{ccc}0&0&\beta^1_s\\0&0&\beta^2_s\\-\beta^1_s&-\beta^2_s&0\end{array}\right)\label{dist}
\end{equation}
where the $\beta^i$'s are independent standard one dimensional brownian motions starting at $0$. The weak convergence of $(B^N_k,Z^N_k, R^N_k)$ to $(B(s), Z(s), R(s))$  is guaranteed by a theorem of Kurtz and Protter \cite{KP} and the identification of $Q_s(\ell_p)=Z(s)e_3$ as a spherical brownian motion is standard exercise in stochastic calculus.\\

\noindent At this point it is useful to recall that each rotation in $SO(3)$ can be identified with an orthonormal frame in $\mathbb R^3$, namely, the frame obtained by rotation of a fixed reference frame. Under this identification the SDE for $Z$ above can be understood as an equation for a moving frame. This SDE is not a stochastic version of the Frenet-Serret frame equations, as one might have thought at the outset, but rather it is a stochastic version of the parallel transport or Bishop frame equations \cite{B}. (Experts in stochastic differential geometry will recognize the stochastic Bishop frame $Z(s)$ as a brownian motion on the orthonormal frame bundle of $S^2$ and $Z(s)e_3$ as the bundle projection onto the base).  Even in the deterministic case, the Bishop frame equations are more general than the Frenet-Serret frame equations in the sense that the former require less regularity of the underlying curve than the latter. This is crucial in the present case, since the Bishop frame parameters that define the Lie algebra process $B(s)$ are a pair of independent one-dimensional brownian motions whereas neither the curvature nor the torsion of $R(s)$ exists in any obvious sense.\\

\noindent It is possible to prove the convergence of the freely rotating chain by adapting the results of Pinsky on isotropic transport processes \cite{Pinsky1}, \cite{Pinsky2}. However, we prefer to develop an approach based on stochastic Bishop frames because it is self-contained and allows us easily to derive the asymptotic behavior of the Kratky-Porod model as the persistence length tends to zero or infinity.


\section{Convergence Theorems}


\subsection{Construction of the Model}
Consider a polymer chain of $N$ segments of length $a$, with each pair of consecutive segments meeting at the same planar angle.  Denote by $\theta$ the angle between two consecutive bonds (see figure $\ref{conformation} $).  Denote the vectors comprising the segments of the polymer by $Q_1,Q_2,\ldots,Q_N$, and let the position of each bead on the polymer, or bond between segments, be denoted by $R_0,R_1,\ldots,R_N$.\\

\noindent Without loss of generality, we can specify the position and orientation of the polymer in $\mathbb{R}^3$ by choosing our coordinates appropriately.  Fix one end of the polymer (the \emph{zeroth bead} or fixed end) at the origin, so that $R_0=0$.  It then follows that
\begin{equation}
R_n=\sum_{i=1}^n Q_i.
\end{equation}
Next, orient the polymer so that the first segment always points straight up, in the positive $z$-direction.  That is,
\begin{equation}
R_1=Q_1=ae_3.
\end{equation}
\noindent Pinning down of one end of the polymer specifies boundary conditions for the model. Since the bond angle $\theta$ is fixed, the bond $Q_n$ can be expressed as a rotation of $Q_{n-1}$ through angle $\theta$ about an axis chosen uniformly at random from the unit vectors in the plane perpendicular to $Q_{n-1}$ . We denote that rotation by $G_n$  and note that in the freely rotating chain, the $G_n$'s are assumed to be mutually independent. Let $Z_n=G_n \, G_{n-1}\cdot ... \cdot G_1$ and note that $Z_n=G_nZ_{n-1}$. Since $Q_n=G_nQ_{n-1}$, it follows that $Q_n=Z_nQ_1$. Thinking of $Z_{n-1}$ as an orthogonal change of basis, we see that the conjugate rotation 

\begin{equation}
H_n=Z_{n-1}^{-1} G_n Z_{n-1}
\end{equation}

\noindent is a rotation of $Q_1$ through an angle $\theta$ about an axis chosen uniformly at random from the unit vectors in the plane perpendicular to $Q_1=ae_3$, and that the $H_n$'s are mutually independent (see figure 2).\\

\begin{figure}[h!]
\includegraphics[scale=0.72]{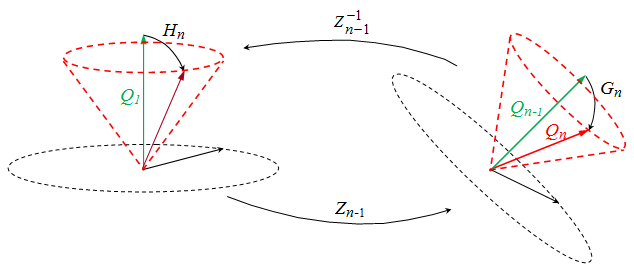}\label{commute}
\caption{Commutative diagram for $H_{n}$ and $G_{n}$}
\end{figure}

\noindent Thus,

\begin{align}
Z_n &= G_n Z_{n-1} \\
 &= Z_{n-1}Z_{n-1}^{-1}G_nZ_{n-1} \\
 &= Z_{n-1}H_n \\
 &= H_1 H_2 \cdot ...\cdot H_n.
\end{align}

\noindent In the sequel, we work with an explicit model of the freely rotating chain by defining $H_n$ to be the rotation through an angle $\theta$ about the axis $(\cos \phi_n, \sin\phi_n, 0)$ where $\phi_1,...,\phi_N$ are i.i.d. random variables uniformly distributed on $[0,2\pi)$. \\

\noindent Because of the product structure of $Z_n$, the freely rotating chain satisfies a simple system of discrete stochastic equations

\begin{align}
R_n  &= \sum_{k=1}^n Z_k (ae_3) \\
Z_n - Z_{n-1} &= Z_{n-1} (H_n-I) \\
R_0 &= 0,\, Z_0=I.
\end{align}

\noindent To analyze the system in the limit $N\to \infty$, we observe that if $N$ is large, then $\theta=\kappa/\sqrt{N}$ is small, hence $H_n$ is close to the identity. Such rotations can be expressed as matrix exponentials of the form

\begin{equation}
H_n = \exp \left(\frac{\kappa}{\sqrt{N}}\, b_n \right)
\end{equation}

\noindent where

\begin{equation}
b_n=\left(\begin{array}{ccc}
0 & 0 & -\sin \phi_n \\
0 & 0 & \cos \phi_n \\
\sin \phi_n & -\cos \phi_n & 0
\end{array}\right)
\end{equation}

\noindent This observation plays a key role in the following result.\\

\begin{thm}
\label{thm (Z^N,B^N) converges to (Z,B)}
Let $\displaystyle a=\frac{L}{N}$, $\displaystyle \theta= \frac{\kappa}{\sqrt{N}}$ and $\displaystyle \ell_p=\frac{2L}{\kappa^2}$. For $0\le s \le L$, define 

$$\displaystyle Z^N_s=Z^N_{[Ns/L]} \:, \text{and}\:  B^N_s=\sum_{k=1}^{[Ns/L]}\left(H_n-I\right).$$

\noindent Then, the process $((Z^N_s,B^N_s);\, 0\le s \le L)$ converges weakly to the process $((Z_s,B_s);\, 0\le s \le L)$, where

\begin{equation}
B_s=\left(\begin{array}{ccc}
0 & 0 & \beta^1_s \\
0 & 0 & \beta^2_s \\
-\beta^1_s & -\beta^2_s & 0 \\
\end{array}\right)
\end{equation}
\noindent and $\beta^1$ and $\beta^2$ are independent standard one dimensional Brownian motions starting at zero. Furthermore, $Z_s$ satisfies the Stratonovich equation

\begin{equation}
\label{dZ = const Z dB}
\partial Z_s=\frac{1}{\sqrt{\ell_p}}Z_s\partial B_s\\
\end{equation}

\noindent with initial condition $Z_0=I$.

\end{thm}

\noindent \textit{Proof}: This theorem is proved in several steps. First, we show that the non-zero entries of $B^N_s$ converge to Brownian motions whose variance at time $s$ is $s/\ell_p$. Second, we use a theorem of Kurtz and Protter \cite{KP} to show that the pair $(B^N_s,Z^N_s;\,0\le s \le L)$ converges weakly to $\displaystyle{\left(\frac{1}{\sqrt{\ell_p}}B_s,Z_s;\,0\le s \le L\right)},$ where $Z_s$ satisfies the Stratonovich equation
$$ \partial Z_s = \frac{1}{\sqrt{\ell_p}}Z_s \partial B_s.$$
\noindent For the first step, note that 
\begin{align} 
\label{B^N}
B^N_s = \sum_{n=1}^{\left[\frac{Ns}{L}\right]}\left(\frac{\kappa}{\sqrt{N}}b_n+\frac{\kappa^2}{2N}b_n^2 + r_n\right)
\end{align}
\noindent where $\displaystyle \|r_n\|\le \frac{e^{\kappa/\sqrt{N}}\kappa^3}{6N^{3/2}}.$\\

\noindent Note that the non-zero entries of $b_n$  are mutually independent random variables with mean 0 and variance $\displaystyle \frac{\kappa^2}{2N}$. Therefore, by Donsker's theorem, the non-zero entries in the first term of $(\ref{B^N})$ converge to independent Brownian motions whose variance at time $s$ is $\displaystyle \frac{s\kappa^2}{2L}=\frac{s}{\ell_p}$. Next, by the Law of Large Numbers, the second term in $(\ref{B^N})$ converges to $\displaystyle \frac{s}{\ell_p}D$, where $D=\text{diag}(-1/2,-1/2,-1)$. Finally, the third term in $(\ref{B^N})$ is $O(N^{-3/2})$ and hence converges to zero. Therefore, the process $\displaystyle (B^N_{\left[\frac{Ns}{L}\right]})$ is tight and converges weakly to $\displaystyle{\left(\frac{1}{\sqrt{\ell_p}}B_s;\, 0\le s \le L\right)}$. Note that $B^N_s$ and $Z^N_s$ are piecewise constant and have jump discontinuities only at integer multiples of $\frac{L}{N}$. Now, if $\displaystyle s=\frac{nL}{N},$ then
\begin{align}
\Delta Z^N_s &= Z^N_s - Z^N_{s^-}\\  
 &= Z^N_n - Z^N_{n-1} \\
 &= Z^N_{n-1}(H^N_n-I) \\
 &= Z^n_{s-}\Delta B^N_s.
\end{align}
\noindent Hence, $Z^N_s$ satisfies the Stratonovich equation
\begin{align}
\partial Z^N_s = Z^N_{s-}\partial B^N_s.
\end{align}
\noindent But by $(\ref{B^N})$, we have $[B^N,B^N]_s=Cs$ for a constant $C$ independent of $N$. \\

\noindent Thus, the sequence $(B^N_s,Z^N_s;\, 0\le s \le L)$ is tight in the Skorohod topology. Consequently, by the theorem of Kurtz and Protter \cite{KP}, this sequence converges weakly to $\displaystyle (\frac{1}{\sqrt{\ell_p}},Z_s;\,0\le s \le L)$ where 
\begin{align}
\partial Z_s = \frac{1}{\sqrt{\ell_p}}Z_{s-}\partial B_s
\end{align}
\noindent which is equivalent to equation $(\ref{dZ = const Z dB})$ since $B_s$ is continuous.

\begin{thm}
\label{thm Q_s is spherical BM}
Let $\displaystyle a=\frac{L}{N}$, $\displaystyle \theta= \frac{\kappa}{\sqrt{N}}$ and $\displaystyle \ell_p=\frac{2L}{\kappa^2}$, and define the process $Q_s=Z_se_3$, $0\le s \le L$. Then, $Q_s(\ell_p)$ is a spherical Brownian motion generated by $\displaystyle \frac{1}{\ell_p}\Delta_{S^2}$. Consequently, the freely jointed chain $(R^N_{[Ns/L]};0\le s \le L)$ converges weakly to the process $\displaystyle R_s=\int_0^s Z_ue_3\;du,\: 0\le s \le L$.
\end{thm}

\noindent \textit{Proof}: The strategy is to use L\'{e}vy's characterization theorem, by which it suffices to show that for any $\mathcal{C}^2$ function $f:S^2\rightarrow \mathbb{R}$, the process

$$f(Q(s))-f(Q(0))-\int_0^s\frac{1}{\ell_p}\Delta_{S^2} f(Q(u)) du $$

\noindent is a martingale. Let $f: S^2\to \mathbb{R}$ be a smooth function and let $F:\mathbb{R}^3\to \mathbb{R}$ be a smooth extension of $f$ such that for some $\varepsilon >0$ and all $\xi \in S^2$, 

\begin{equation}
F(r\xi ) = f(\xi), \hspace{.3cm} 1-\varepsilon < r < 1+\varepsilon.
\end{equation}

\noindent Then, the radial derivative of $F$ vanishes, and so the expression of the Laplacian in $\mathbb{R}^3$ in polar coordinates simplifies to 

\begin{equation}
\Delta_{\mathbb{R}^3}F(\xi)=\Delta_{S^2}F(\xi)=\Delta_{S^2}f(\xi)
\end{equation}

\noindent for all $\xi \in S^2$.\\

\noindent Let $X_s=Z_se_1$, $Y_s=Z_se_2$ and $Q_s=Z_se_3$. Then the matrix equation $\partial Z_s=Z_s\partial B_s$ is equivalent to the system

\begin{align}
\partial X_s &= -\frac{1}{\sqrt{\ell_p}}\, Q_s \partial \beta^1_s \\
\partial Y_s &= -\frac{1}{\sqrt{\ell_p}}Q_s\partial \beta^2_s \\
\partial Q_s &= \frac{1}{\sqrt{\ell_p}}\big(X_s \partial \beta^1_s + Y_s \partial \beta^2_s\big).
\end{align}

\noindent Applying It\^o's formula to $Q_s$ thought of as a process in $\mathbb{R}^3$ yields

\begin{align}
f(Q_s)-f(e_3) &= F(Q_s)-F(e_3) \\
 &= \frac{1}{\sqrt{\ell_p}} \int_0^s \nabla F(Q_s)\cdot \partial Q_s \\
 &= \frac{1}{\sqrt{\ell_p}}\int_0^s \nabla F(Q_u)\cdot dQ_u + \frac{1}{2\sqrt{\ell_p}}\int_0^s \text{Hess}\, F(Q_u):d\langle Q, Q \rangle_u 
\end{align}

\noindent The first term above is a martingale. To evaluate the integrand of the second term, note that 

\begin{equation}
d\langle Q, Q \rangle_u = (X_u \otimes X_u)du 
\end{equation}
\noindent hence 

\begin{align}
\text{Hess}\, F(Q_u):d\langle Q, Q\rangle_u  
=\frac{\partial^2}{\partial X^2_u}F(Q_u) + \frac{\partial^2}{\partial Y^2_u}F(Q_u)
\end{align}

\noindent where the directional derivatives are computed $\omega$ by $\omega$. Since the unit vector $Q_u$ points almost surely in the radial direction, we have

\begin{align}
\frac{\partial{\partial^2}}{\partial Q_u^2} F(Q_u)
&= \frac{\partial^2}{\partial r^2} F(Q_u) = 0 \hspace{.2cm} \text{a.s.}
\end{align}

\noindent and hence

\begin{align}
\text{Hess}\, F(Q_u):d\langle Q,Q\rangle_u &= \Big( \frac{\partial^2}{\partial X^2_u} + \frac{\partial^2}{\partial Y^2_u} + \frac{\partial^2}{\partial Q^2_u} \Big) F(Q_u) \\
&= \Delta_{\mathbb{R}^3}\hspace{.2cm}F(Q_u) \text{ a.s.} \\
&= \Delta_{S^2}f(Q_u).
\end{align}

\noindent Here we have used the fact that

\begin{equation}
\Delta_{\mathbb{R}^3}F(x)=\Big( \frac{\partial^2}{\partial u^2_u} + \frac{\partial^2}{\partial v^2_u} + \frac{\partial^2}{\partial w^2_u} \Big) F(x)
\end{equation}
 
\noindent for any orthogonal frame $(u,v,w)$ based at $x\in \mathbb{R}^3$.

\section{Asymptotic Behavior in the Persistence Length}

One advantage of the representation of the tangent process $Q_s=\displaystyle \frac{dR_s}{ds}=Z_se_3$ in terms of the $SO(3)$-valued process $Z$ with $\partial Z_s=\frac{1}{\sqrt{\ell_p}}Z_s\partial B_s$ is that we can analyze the behavior of the Kratky-Porod model when the persistence length $\ell_p$ tends to $0$ or infinity. We first need a lemma about the correlation between vectors $Q_s$ and $Q_t$ at different ``times" (i.e. arc lengths) $s$ and $t$.

\begin{lemma}
\label{lemma E[Qs*Qt]=exp(t-s)}
Let $Q_s=Z_s e_3$ be the process as above. Then $\displaystyle E[Q_s\cdot Q_t]=e^{-\frac{2|t-s|}{\ell_p}}.$ 
\end{lemma}

\noindent \textit{Proof}: Suppose $s<t$. Since $Q_s\cdot Q_t$ is rotationally invariant, it follows that for any $\xi\in S^2$, and rotation $O\in SO(3)$ with $O\xi=e_3$,

\begin{align}
E_{\xi} [Q_0\cdot Q_s]&= E_{\xi}[OQ_0\cdot OQ_s]\\
&= E_{e_3}[e_3\cdot Q_s]
\end{align}

\noindent By the Markov property, $
E_{e_3}[Q_s\cdot Q_t] = E_{e_3}[E_{Q_s}[Q_0\cdot Q_{t-s}\,]]= E_{e_3}[Q_0\cdot Q_{t-s}]$.

\noindent Therefore, to prove the lemma, it suffices to show $\displaystyle E[Q_s\cdot e_3]=e^{-\frac{2s}{\ell_p}} $, for all $s\in[0,L]$.\\

\noindent To this end, define $f:\mathbb{R}^3\rightarrow \mathbb{R}$ by $f(P)=P\cdot e_3$. Alternatively, $f$ has the form $f(P)=r\cos\phi$ in polar coordinates of the point $P=(r, \theta,\phi)$.  Then

\begin{align}
\label{eqn Delta_R^3 to Delta_S^2}
\Delta_{\mathbb{R}^3}f(P) &=
\frac{\partial^2}{\partial r^2}f(P)+\frac{2}{r}\,\frac{\partial}{\partial r}f(P) + \frac{1}{r^2}\,\Delta_{S^2}f(P)\\
&= \frac{2}{r}\cos\phi + \frac{1}{r^2}\,\Delta_{S^2}f(P).
\end{align}

\noindent On the other hand,
 $\Delta_{\mathbb{R}^3}f\equiv 0$, as $f$ is a linear function of cartesian coordinates of its argument. Thus,

\begin{align}
\Delta_{S^2}f(P)= - 2\, r \cos\phi = -2 f(P).
\end{align}

\noindent Now, applying It\^{o}'s formula for $Q_s$, we find that the process

\begin{align}
M_s &= f(Q_s(\ell_p)) -  f(Q_0(\ell_p)) - \int_0^s\frac{1}{\ell_p}\Delta_{S^2}f(Q_u(\ell_p))du\\
&= \label{eqn 350} f(Q_s(\ell_p)) -1  + \frac{2}{\ell_p} \int_0^s f(Q_u(\ell_p))du
\end{align}

\noindent is a martingale with mean 0. Taking the expectation in $(\ref{eqn 350})$, we conclude that the non-random function $x_s=E[f(Q_s(\ell_p))]$ is the solution of the initial value problem

\begin{align}
\dot{x}_s &= -\frac{2}{\ell_p}x_s,\hspace{.2cm} x_0 = 1;
\end{align}


\noindent hence, $E[Q_s\cdot e_3] = x_s = e^{-\frac{2s}{\ell_p}}$.

\vspace{.5cm}

\begin{lemma}
\label{lemma E[|R_t|^2] is prop. to s}
$\displaystyle E\left[|R_t|^2\right] = \ell_p\, t - \frac{\ell_p^2}{2}\,\left(1-e^{-\frac{2t}{\ell_p}}\right).$
\end{lemma}

\noindent \textit{Proof}: 
\begin{align}
E\left[|R_t|^2\right] &= E\left[ \int_0^tQ_u\,du\, \cdot\, \int_0^t Q_v\, dv\right] \\
&= 2\, \int_0^t \int_0^v E\left[ Q_u\cdot Q_v \right]\, du\, dv\\
&= 2\, \int_0^t \int_0^v e^{\frac{-2(v-u)}{\ell_p}}du\,dv\\
&= \ell_p t - \frac{\ell_p^2}{2}\,\left(1-e^{-\frac{2t}{\ell_p}}\right).
\end{align}

\vspace{.5cm}
\begin{thm}
\label{Thm 2 repeated (asympt. behavior of l_p)}
The Kratky-Porod model with fixed length $L$ exhibits the hard rod/random coil transition, in the sense that,

A. (Hard Rod) As $\ell_p\to\infty$, $(R_s; \, 0\le s \le L)\, \to \, (se_3;\, 0\le s \le L)$ in probability, and $(\sqrt{\ell_p}(R_s-se_3);\, 0\le s \le L)$ converges weakly to a Gaussian process $(W_s;\, 0 \le s \le L)$ with components $W^3_s\equiv 0$, and for $i=1,2$,

$$
W^i_s=\int_0^s\beta^i_udu,~ 0\leq s\leq L  
$$

\noindent where the $\beta^i$'s are independent, standard, one dimensional brownian motions starting at zero.

B. (Random Coil) As $\ell_p\to 0$,  $(\sqrt{3/\ell_p} R_s;\, 0\le s \le L)$ converges weakly to $(W_s;\, 0 \le s \le L)$, where $W$ is a standard 3-dimensional Brownian motion.

\end{thm}

\vspace{1cm}

\noindent \textit{Proof of A.} Since 
\begin{align}
 R_s - se_3 &= \int_0^s (Z_u-I)e_3\, du,
\end{align}
it is enough to show that $\displaystyle \int_0^{\cdot}(Z_u - I)\,du\to 0 $ in probability as $\ell_p\to \infty$. To this end, for a matrix-valued process $M=[m^{ij}]_{i,j=1}^3$ we define a norm $||M ||$ pathwise by the formula 
\begin{equation}\| M\|^2 = \displaystyle \sum_{i,j=1}^3 \sup_{0\le s \le L}|m^{ij}_s|^2
\end{equation} 

\noindent and we will prove that the random variable
$\displaystyle \left\|\int_0^{\cdot}(Z_u-I)\,du\right\|$ converges to 0 in probability. Now, it is not hard to see that every integrable matrix-valued process $M$ satisfies
\begin{eqnarray}
\label{norm(integral)<L norm(integrand)}
\left\|\int_0^{\cdot}M_u\,du\right\| \le L \|M\|.
\end{eqnarray}
 
\noindent Thus, it suffices to show $\| Z-I\|\to 0$ in probability as $\ell_p\to \infty$.\\
For $\varepsilon>0$, we have

\begin{eqnarray}
P\{ \, \|Z-I\|\ge \varepsilon \} &=& P\left\{ \, \left\|\int_0^{\cdot} Z_u\partial B_u\, \right\|\ge \varepsilon \sqrt{\ell_p} \right\}\nonumber\\
&\le & \label{ineq Chebishev for integral ZdB} \frac{1}{\varepsilon^2\ell_p}E\left[\left\| \int_0^{\cdot}Z_u\partial B_u \, \right\|^2\right]
\end{eqnarray}
and so it suffices to show that the expected value on the right hand side of ($\ref{ineq Chebishev for integral ZdB}$) is bounded. Note that the matrix $\displaystyle \int_0^s Z_u\, \partial B_u$ has entries of the form $\displaystyle \int_0^s z^{ij}_u\,\partial \beta^k_u$, for all $0\le i,j\le 3$ and $k=1,2$, or perhaps a sum of two such integrals. Assuming 
\begin{equation} C = \max_{i,j,k}E\left[\sup_{0\le s \le L}\left|\int_0^s z^{ij}_u\, \partial\beta^k_u\, \right|^2\right] \, < \infty,
\end{equation}

\noindent then 
\begin{equation}
E\, \left\|\int_0^{\cdot}Z_u\,\partial B_u \, \right\|\le 36C,
\end{equation}
\noindent as one can see by squaring out the terms of the sum in the definition of the norm $\|M\|$ and estimating the cross term by the Cauchy-Schwarz inequality.\\

\noindent Thus, to prove that the expected value on the right hand side of $(\ref{ineq Chebishev for integral ZdB})$ is bounded, it suffices to show that 

\begin{equation}
\label{Esup_square is bdd}
 E\left[\sup_{0\le s \le L}\left|  \int_0^s z_u\partial \beta_u\, \right|^2 \right]<\infty
\end{equation} 
\noindent (where we have dropped the superscripts for ease of notation). The conversion of Stratonovich to \^Ito integrals gives us:

\begin{equation}
\label{Strato to Ito}
\int_0^s z_u\partial \beta_u = \int_0^s z_u\,d\beta_u +\frac{1}{2}\langle z,\beta \rangle_s.
\end{equation}

\noindent Observe that \\

\noindent $\displaystyle
\sup_{0\le s \le L}\left| \int_0^s z_u\,d\beta_u + \frac{1}{2}\langle z,\beta \rangle_s\right|^2 $
\begin{align}
\le \sup_{0\le s \le L}\left| \int_0^s z_u\,d\beta_u\right|^2 &+  \sup_{0\le s \le L}\left| \int_0^s z_u\,d\beta_u\right|\cdot \sup_{0\le s \le L}\big|\langle z,\beta \rangle_s\big|  \nonumber \\
&+ \frac{1}{4}\,\sup_{0\le s \le L}\big|\langle z,\beta \rangle_s \big|^2.
\end{align}

\noindent Also, by Jensen's inequality,

\begin{equation}
E\left[\sup_{0\le s \le L}\left| \int_0^s z_u\,d\beta_u \right|\right] \le E\left[\sup_{0\le s \le L}\left| \int_0^s z_u\,d\beta_u \right|^2\right]^{1/2}.
\end{equation}\\

\noindent Therefore, to prove $(\ref{Esup_square is bdd})$, it suffices to show 
\begin{equation}
\label{Esup_square of Ito integ is bdd} E\left[\sup_{0\le s \le L}\left| \int_0^s z_u\,d\beta_u\right|^2\right]<\infty
 \end{equation}
and that for some constant $K\ge 0$, almost surely,
\begin{equation}
\label{sup quad var is bdd} \sup_{0\le s \le L}\big| \langle z,\beta\rangle_s \big| \le K.
 \end{equation}
By Doob's inequality, and using the fact that $|z_u|\le 1$ almost surely,

\begin{eqnarray}
E\left[\sup_{0\le s \le L} \left| \int_0^s z_u\,d\beta_u\right|^2\right] &\le & 4E\left[\int_0^L z_u^2\,du\right] \\
&\le & 4 L
\end{eqnarray}

\noindent which shows $(\ref{Esup_square of Ito integ is bdd})$. To prove $(\ref{sup quad var is bdd})$, we note that $\langle z, \beta\rangle_s$ can be written as the sum of at most two terms of the form $\displaystyle \ell_p^{-1/2} \left\langle \int_0^{\cdot} \tilde{z}_u\partial\beta^i_u,\beta \right\rangle$, with $\tilde{z}$ being an entry of $Z$. Using $d\langle\beta^i,\beta^j\rangle_u = \delta_{ij}\,du$, it follows that $\langle z, \beta \rangle_s$ is equal to only one term of the form $\displaystyle \ell_p^{-1/2} \int_0^s \tilde{z}_u du$. Thus,

\begin{align}
\sup_{0\le s \le L}\big| \langle z,\beta \rangle_s  \big| &\le  \frac{1}{\sqrt{\ell_p}}\,\sup_{0\le s \le L}\int_0^s \,\max_{ i,j}|z^{ij}_u|\, du  \\
& \label{quad cov. is bdd} \le\frac{L}{\sqrt{\ell_p}}
\end{align}

\noindent which completes the proof of the first part of $A$.\\

\noindent Recalling the definition of the Gaussian proces $W_s, 0\le s \le L$, let us observe that

\begin{align}
\sqrt{\ell_p}(R_s-se_3)-W_s &= \sqrt{\ell_p}\int_0^s (Z_u-I)e_3\,du - \int_0^s B_u e_3\,du \\
\label{integ of difference}&= \int_0^s \left(\sqrt{\ell_p}(Z_u-I)-B_u \right)e_3\,du.
\end{align}
Let us show that the norm of the expression in  $(\ref{integ of difference})$ converges to 0 in probability, as $\ell_p\to \infty$.
However, by  $(\ref{norm(integral)<L norm(integrand)} )$, it is enough to show that

\begin{align}
\big\| \sqrt{\ell_p}(Z-I)-B\big\| \to0\hspace{.2cm}
\end{align}

\noindent in probability. 

\noindent Now, iterating the equation $\partial Z_s=\displaystyle \ell_p^{-1/2}Z_s\partial B_s$ twice, it follows that for all $s\in[0,L]$,
\begin{eqnarray}
Z_s-I &=& \frac{1}{\sqrt{\ell_p}} \int_0^s\left( I + \frac{1}{\sqrt{\ell_p}} \int_0^u Z_v\,\partial B_v \right)\,\partial B_u \\
&=& \frac{1}{\sqrt{\ell_p}}B_s + \frac{1}{\ell_p}\int_0^s\int_0^u Z_v\,\partial B_v\, \partial B_u 
\end{eqnarray}

\noindent which yields

\begin{eqnarray}
\sqrt{\ell_p}(Z_s-I)-B_s=\frac{1}{\sqrt{\ell_p}}\int_0^s\int_0^uZ_v\,\partial B_v\, \partial B_u.
\end{eqnarray}

\noindent Therefore, for all $\varepsilon>0$,

\begin{eqnarray}
P\left\{ \|\sqrt{\ell_p}\big(Z-I\big)-B \| \ge \varepsilon \right\} \le \frac{1}{\varepsilon^2\ell_p}E\left[\sum_{i,j=1}^3 \sup_{0\le s \le L}|x^{ij}_s|^2 \right]
\end{eqnarray}

\noindent where $[x^{ij}_s]=X_s=\displaystyle \int_0^s \int_0^u Z_v\, \partial B_v\, \partial B_u $. As in the proof of the first part of $A$, we observe that each entry $x^{ij}_s$ of the matrix $X_s$ is the sum of at most two terms of the form

\begin{equation}
\label{form of x}
\int_0^s \int_0^u z_v
\,\partial
 \beta^k_v
 \,\partial \beta^l_u
\end{equation}

\noindent where $z_v$  is an entry of $Z_v$. To complete the proof, it suffices to show that

\begin{equation}
\label{E sup double integ is bdd}
E\left[\sup_{0\le s \le L}\left|\int_0^s \int_0^uz_v\,\partial \beta^k_v\,\partial \beta^l_u \right|^2\right] < \infty.
\end{equation} 
This is obvious if $x^{ij}$ is equal to just one term of the form $(\ref{form of x})$. If, however, $x^{ij}$ is a sum of two such terms, the sufficency of  $(\ref{E sup double integ is bdd})$ becomes clear after expanding $ E[\sup|x_s^{ij}|^2]$ and using the Cauchy-Schwarz inequality.\\

\noindent Now, to see $(\ref{E sup double integ is bdd})$, we convert from Stratonovich to It\^{o} integrals and obtain

\begin{align}
\displaystyle \int_0^s \int_0^u z_v\,\partial\beta^k_v\,\partial\beta^l_u 
&= \int_0^s\left(\int_0^u \,z_v\,d\beta^k_v+\frac{1}{2}\langle z,\beta^k\rangle_s \right)\,\partial\beta^l_u\hspace{2cm}\\
&= \int_0^s\int_0^u z_v\,d\beta^k_v\,d\beta^l_u + \frac{1}{2}\left\langle \int_0^{\cdot}z_v\,d\beta^k_v,\, \beta^l \,\right\rangle_s\nonumber\\
& \label{eqn expansion of double integral} \hspace{1cm} + \frac{1}{2}\int_0^s\langle z,\beta^k\rangle_u\, d\beta^l_u + \frac{1}{4}\left\langle \langle z,\beta \rangle , \, \beta^k \,\right\rangle_s.
\end{align}

\noindent In fact, $\big\langle \langle z,\beta\rangle , \, \beta^k\rangle_s = 0$, since $\langle z, \beta\rangle$ is of finite variation. Also, 

\begin{align}
\sup_{0\le s \le L}\left|\,\left\langle \int_0^{\cdot}z_v\,d\beta^k_v,\, \beta^l \,\right\rangle_s  \right| &= \delta_{kl}\cdot \sup_{0\le s \le L}\left|\,\int_0^s z_v\,dv \,  \right| \\
&\label{sup quad integ is bdd}\le L
\end{align}

\noindent almost surely. Furthermore, by Doob's inequality and $(\ref{quad cov. is bdd})$,

\begin{align}
\displaystyle E\left[\sup_{0\le s \le L}\left|\int_0^s\langle z,\beta^k\rangle_u\, d\beta^l_u\right|^2 \right]
&\le  4\, E\left[\int_0^L \big| \langle z,\, \beta^k\rangle_u \big|^2\,du\right] \\
\label{integ of quad var is bdd}&\le \frac{4L^3}{\ell_p}.
\end{align}

\noindent Therefore, squaring $(\ref{eqn expansion of double integral})$, and using $(\ref{quad cov. is bdd})$, $(\ref{sup quad integ is bdd})$ and $(\ref{integ of quad var is bdd})$, we see that to prove $(\ref{E sup double integ is bdd})$ for Stratonovich integrals, it suffices to prove the boundedness of the corresponding quantity for \^Ito integrals, that is,

\begin{equation}
E\left[\sup_{0\le s \le L}\left|\int_0^s \int_0^uz_v \, d\beta^k_v \, d\beta^l_u \right|^2\right] < \infty.
\end{equation}

\noindent But, by Doob's inequality,
\begin{align} 
E\left[\sup_{0\le s \le L}\left|\int_0^s\int_0^u z_v\,d\beta^k_v\,d\beta^l\right|^2\right]
 &\le 4 \cdot E\left[\int_0^L  \left| \int_0^u z_vd\beta^k_v \right|^2 du\right]\\
&\le  4 L \cdot E\left[ \sup_{0\le u \le L}\left| \int_0^u z_v\, d\beta^k_v \right|^2 \right]\\
&\le  16 L \cdot E\left[\int_0^L |z_u|^2 \, du \right]\\
&\le  16 L^2.
\end{align}

\noindent This completes the proof of part $A$.\\

\noindent \textit{Proof of B}. We would like to prove a functional central limit theorem for $(R_s;\, 0\le s \le L)$ by viewing $R_s=(R^1_s,R^2_s,R^3_s)$ as a vector of additive functionals of the underlying spherical Brownian motion:

\begin{align}
R^i_s(\ell_p)&=\int_0^sQ^i_u(\ell_p)\,du \\
&= \int_0^s f^i(Q_u(\ell_p))\, du
\end{align}

\noindent where $f^i: S^2\rightarrow \mathbb{R}$ is the $i^{\mathrm{th}}$ coordinate function of $Q$, i.e. $f^i(Q)=Q^i$, for $i=1,2,3$.\\

\noindent Let $g^i:S^2\rightarrow \mathbb{R}$ be a solution of the Poisson equation $\Delta_{S^2}g^i = f^i$, and let $\sigma$ be the normalized Lebesgue measure on the unit sphere. By symmetry, $\displaystyle \int_{S^2}f^i(Q)\sigma(dQ)=0$, and since $f^i$ is smooth, there exists a smooth solution $g^i$, unique up to an additive constant. By a change of variable, the process $\displaystyle (Q_s(\ell_p);\, s\ge 0)$ is equal in distribution to $ (Q_{\frac{s}{\ell_p}}(1);\, s\ge 0)$, as well as

\begin{equation}
\label{R_s(lp)=lp R_s/lp(1)}
 (R_s(\ell_p);\, s\ge 0)\stackrel{\mathcal{D}}{=} (\ell_p R_{\frac{s}{\ell_p}}(1);\, s\ge 0).
 \end{equation}

\noindent Using It\^{o}'s formula and $\Delta_{S^2}g^i=f^i$, we find that
\begin{align}
 g^i(Q_s(\ell_p)) - g^i(Q_0(\ell_p)) -  \int_0^s \frac{1}{\ell_p} f^i(Q_u(\ell_p))\,du
\end{align}
\noindent is a martingale. So,
\begin{align}
\label{R is a mtgle plus negligible term}
\frac{1}{\ell_p}R^i_s(\ell_p) = \int_0^s \frac{1}{\ell_p} f^i(Q_u(\ell_p))\,du=g^i(Q_s(\ell_p))-g^i(Q_0(\ell_p))+M^i_s(\ell_p)
\end{align}
\noindent for some martingale $M^i_s(\ell_p)$, $1\le i \le 3$. In particular, for any $t\ge 0$,

\begin{align}
\label{R(1) is a mtgle plus negligible term}
R^i_t(1) =g^i(Q_t(1))-g^i(Q_0(1))+M^i_t(1).
\end{align}

\noindent Comparing $(\ref{R(1) is a mtgle plus negligible term})$ for $t=s/\ell_p$ with $(\ref{R is a mtgle plus negligible term})$, and using the fact that $Q_s(\ell_p)\stackrel{\mathcal{D}}{=}\displaystyle Q_{\frac{s}{\ell_p}}(1)$ along with $(\ref{R_s(lp)=lp R_s/lp(1)})$, it follows that $M^i_s(\ell_p)\stackrel{\mathcal{D}}{=}\displaystyle M^i_{\frac{s}{\ell_p}}(1)$. This yields

\begin{align}
\frac{1}{\sqrt{\ell_p}}\, R^i_{s}(\ell_p) =\sqrt{\ell_p}\, \Big[ g^i(Q_{\frac{s}{\ell_p}}(1))-g^i(Q_0(1))\Big] +\sqrt{\ell_p}\, M^i_{\frac{s}{\ell_p}}(1).
\end{align}

\noindent In fact, $M^i_s(1)$, which we denote for short by $M^i_s,\, 1\le i \le 3$, is a system of continuous martingales with predictable quadratic covariation

\begin{equation}
\langle M^i, M^j \rangle_{t}=  \int_0^{t}\Gamma (g^i,g^j)(Q_s(1))\, ds
\end{equation}

\noindent where $\Gamma$ is the carr\'{e} du champ operator:

\begin{equation}
\Gamma(\varphi,\psi)=\frac{1}{2}\Big[\Delta_{S^2}(\varphi\psi)(Q)-\varphi(Q)\Delta_{S^2}\psi(Q)-\psi(Q)\Delta_{S^2}\varphi(Q)\Big].
\end{equation}

\noindent Since $g^i$ is bounded, the functional central limit theorem for the stochastic process 
 $1/\sqrt{\ell_p}R^i_s(\ell_p)=\sqrt{\ell_p}R^i_{\frac{s}{\ell_p}}(1)=\sqrt{\ell_p}R^i_{\frac{s}{\ell_p}}$ as $\ell_p\rightarrow 0$ is equivalent to the corresponding result for the martingale $\displaystyle \sqrt{\ell_p}\,M^i_{\frac{s}{\ell_p}}$. According to the functional central limit theorem for martingales (e.g. \cite{EK}, \cite{W}), this amounts to showing that 

\begin{equation}
\lim_{\lambda\rightarrow \infty}\frac{1}{\lambda}\langle M^i, M^j\rangle_{\lambda t}=c_{ij}t
\end{equation}

\noindent for some non-negative definite matrix $C=[c_{ij}]$. By the ergodic theorem for $Q$, 

\begin{align}
c^{ij}=\lim_{t\rightarrow \infty}\frac{1}{t} \langle M^i,M^j\rangle_t &= \lim_{t\rightarrow \infty}\frac{1}{t}\int_0^{t}\Gamma(g^i,g^j)(Q_s)\,ds \\
&= \int_{S^2}\Gamma(g^i,g^j)(Q)\sigma(dQ)
\end{align}

\noindent almost surely and in $L^1$. Rather than calculating $\Gamma(g^i,g^j)$ explicitly, we can identify $C$ by a symmetry argument. We will show that $OCO^T=C$ for any rotation $O\in SO(3)$, which implies that $C$ is a constant multiple of the identity matrix.\\

\noindent To see this, let $O\in SO(3)$ and define a process $\tilde{R}$ as 
\begin{align}
\tilde{R}_s = OR_s &= \int_0^s OQ_u(1)\, du\\
&=\int_0^s\tilde{Q}_u(1)\,du
\end{align}
where $\tilde{Q}=OQ$ is a spherical brownian motion  generated by $\Delta_{S^2}$, starting at $Oe_3$. Because $\tilde{R}$ is a Kratky-Porod model, $\tilde{R}$ differs from a martingale, say, $\tilde{M}$ by a bounded term, due to equation $(\ref{R is a mtgle plus negligible term})$. Hence,

\begin{align}
\lim_{t\to \infty}\frac{1}{t}E[\tilde{R}_t \tilde{R}_t^T] &= \lim_{t\to \infty}\frac{1}{t}E[\tilde{M}_t \tilde{M}_t^T]\\
&= \lim_{t\to \infty}\frac{1}{t}E \langle \tilde{M}, \tilde{M}\rangle_t \\
&= \lim_{t\to \infty}\frac{1}{t}E\langle M, M\rangle_t = C. 
\end{align}

\noindent On the other hand,

\begin{align}
C = \lim_{t\to \infty}\frac{1}{t}E[\tilde{R}_t \tilde{R}_t^T] &= \lim_{t\to \infty}\frac{1}{t}\,O E[R R^T] O^T \\
&= \lim_{t\to \infty}\frac{1}{t}\,O E\langle M, M\rangle_t O^T \\
&= OCO^T.
\end{align}

\noindent So, $C=cI$ for some constant $c$. Now, using Lemma $\ref{lemma E[|R_t|^2] is prop. to s}$, 

\begin{align}
\lim_{t\rightarrow \infty}\frac{1}{t} E\big[ |\tilde{R}_t |^2 \big] = \lim_{t\rightarrow \infty} \frac{1}{t}E\big[|R_t|^2\big]=1.
\end{align}
Now, 
\begin{align}
\lim_{t\rightarrow \infty }\frac{1}{t}E\big[|R_t|^2\big]=\text{tr}\,C=3c
\end{align}

\noindent and so $c=1/3$. \\
This completes the proof of Theorem $\ref{Thm 2 repeated (asympt. behavior of l_p)}$.

\section{Conclusion}

The Kratky-Porod model for semi-flexible polymers in dilute solution is a continuous model that proves to be the continuum limit of the discrete freely rotating chain model, under an appropriate scaling.  While this result had been known for over sixty years in the field of chemical physics, our Theorem 4 provides the first known rigorous, probabilistic proof of the correspondence between the two models.  This is shown by means of an intermediate result (Theorem 3) that relates the tangent vector along the polymer to a sequence of matrix-valued processes, defined by a stochastic differential equation.  This sequence is shown to converge by a theorem of Kurtz and Protter, and the limiting vector-valued process is found to be a spherical Brownian motion that matches the continuous Kratky-Porod model.\\

\noindent Moreover, the model is proven to converge to the expected results in two extreme cases: for long persistence lengths (large values of $\ell_p$), the adjacent segments are highly correlated, and so the polymer approximates a straight rod; while for short persistence lengths (small values of $\ell_p$), the adjacent segments are nearly independent, and the polymer approaches the Rouse or freely jointed chain model, in which the tangent vector approaches a three-dimensional Brownian motion.  These limiting cases are shown in Theorem 7, which provides the additional result that in the former case, the deviation from the straight rod is itself a Brownian motion in the plane.  This result places the Kratky-Porod model on a continuum with the straight rod and the Rouse model as the stiffness parameter, or persistence length, varies.\\

\noindent With these results firmly in place, we can expand the Kratky-Porod model into other realms.  One such case, in which a constant external force is applied to the polymer, introduces other limiting regimes in which the shape and orientation of the polymer depend both on the persistence length and relative strength and direction of the force.  Another case is the Rotational Isomeric State approximation model, in which conformation of the polymer is determined by finite and fixed number of possible states in each step of the conformation.\\

\noindent These models, therefore, unite the fields of probability theory and chemical physics, using the methods and tools of the former field to demonstrate a property in the latter.  At long last, the relationship between the freely rotating and Kratky-Porod polymer models has been proven.

\end{document}